\documentclass[12pt]{amsart}
\usepackage{amsmath, amssymb, latexsym, amsthm, color} 

\long\def\forget#1\forgotten{}
\forget
\forgotten

\newcommand{\nc}{\newcommand}
\nc{\my}[1]{\textcolor{red}{\sf #1}}
\nc{\op}{\operatorname}\nc{\FP}{\op{FP}}\nc{\FPhat}{\widehat{\op{FP}}}
\nc{\Impl}{\Rightarrow}
\nc{\meet}{\wedge}
\nc{\N}{\mathbb{N}}
\nc{\Z}{\mathbb{Z}}
\nc{\cF}{\mathcal{F}}
\nc{\set}[2]{\{#1 \,:\, #2\}}
\nc{\thusfar}{\marginpar{***}{\sf finished thus far}}
\nc{\DL}{\operatorname{D}_\mathrm{L}}
\nc{\DR}{\operatorname{D}_\mathrm{R}}
\nc{\NDL}{\operatorname{ND}_\mathrm{L}}
\nc{\NDR}{\operatorname{ND}_\mathrm{R}}
\newtheorem{thm}{Theorem}[section]
\nc{\bthm}{\begin{thm}} \nc{\ethm}{\end{thm}}
\newtheorem{prop}[thm]{Proposition}
\nc{\bprp}{\begin{prop}} \nc{\eprp}{\end{prop}}
\newtheorem{fact}[thm]{Fact}
\nc{\bfct}{\begin{fact}} \nc{\efct}{\end{fact}}
\newtheorem{prob}[thm]{Problem}
\nc{\bprb}{\begin{prob}} \nc{\eprb}{\end{prob}}
\newtheorem{lem}[thm]{Lemma}
\nc{\blem}{\begin{lem}} \nc{\elem}{\end{lem}}
\newtheorem{claim}[thm]{Claim}
\nc{\bclm}{\begin{claim}} \nc{\eclm}{\end{claim}}
\newtheorem{cor}[thm]{Corollary}
\nc{\bcor}{\begin{cor}} \nc{\ecor}{\end{cor}}
\newtheorem{conj}[thm]{Conjecture}
\nc{\bcnj}{\begin{conj}} \nc{\ecnj}{\end{conj}}
\theoremstyle{definition}
\newtheorem{defn}[thm]{Definition}
\nc{\bdfn}{\begin{defn}} \nc{\edfn}{\end{defn}}
\newtheorem{cnv}[thm]{Convention}
\nc{\bcnv}{\begin{cnv}} \nc{\ecnv}{\end{cnv}}
\theoremstyle{remark}
\newtheorem{rem}[thm]{Remark}
\nc{\brem}{\begin{rem}} \nc{\erem}{\end{rem}}
\newtheorem{exam}[thm]{Example}
\nc{\bexm}{\begin{exam}} \nc{\eexm}{\end{exam}}
\nc{\bpf}{\begin{proof}} \nc{\epf}{\end{proof}}
\nc{\be}{\begin{enumerate}}
\nc{\ee}{\end{enumerate}}
\nc{\bi}{\begin{itemize}}
\nc{\ei}{\end{itemize}}
\nc{\itm}{\item}
\nc{\sm}{\setminus}
\nc{\sub}{\subseteq}
\nc{\inv}{^{-1}}
\nc{\la}{\langle}
\nc{\ra}{\rangle}

\nc{\ed}{

\end{document}}

\begin{document}

\title[Hindman's Coloring Theorem in semigroups]{Hindman's Coloring Theorem in arbitrary semigroups}

\author{Gili Golan}
\email{gili.golan@math.biu.ac.il}

\author{Boaz Tsaban}
\email{tsaban@math.biu.ac.il}
\urladdr{http://www.cs.biu.ac.il/\~{}tsaban}

\address{Department of Mathematics, Bar-Ilan University, 5290002 Ramat Gan, Israel}



\begin{abstract}
Hindman's Theorem asserts that, for each finite coloring of the natural numbers,
there are distinct natural numbers $a_1,a_2,\dots$ such that all of the
sums $a_{i_1}+a_{i_2}+\dots+a_{i_m}$ ($m\ge 1$, $i_1<i_2<\dots<i_m$) have the same color.

The celebrated Galvin--Glazer proof of Hindman's Theorem and a classification
of semigroups due to Shevrin, imply together that, for each finite coloring of each infinite semigroup $S$, there
are distinct elements $a_1,a_2,\dots$ of $S$ such that all but finitely many of the products
$a_{i_1}a_{i_2}\cdots a_{i_m}$ ($m\ge 1$, $i_1<i_2<\dots<i_m$) have the same color.

Using these methods, we characterize the semigroups $S$ such that,
for each finite coloring of $S$, there is an infinite \emph{subsemigroup} $T$ of $S$, such that
all but finitely many members of $T$ have the same color.

Our characterization connects our study to a classical problem of Milliken,
Burnside groups and Tarski Monsters.
We also present an application of Ramsey's graph-coloring theorem to Shevrin's theory.
\end{abstract}

\maketitle

\tableofcontents

\section{The Galvin--Glazer--Hindman Theorem}

A \emph{finite coloring} of a set $A$ is an assignment of one color to each element of $A$, where the set of possible colors is finite.
In 1974, Hindman proved the following theorem, extending profoundly a result of Hilbert.

\bthm[Hindman \cite{hindman}]\label{hind}
For each finite coloring of $\N$, there are $a_1,a_2,\dots\in\N$ such that all sums $a_{i_1}+a_{i_2}+\dots+a_{i_m}$
($m\ge 1$, $i_1<i_2<\dots<i_m$) have the same color.
\ethm

In Hindman's Theorem \ref{hind}, we may request that the elements $a_1,a_2,\dots$ are \emph{distinct},
by moving, if needed, to appropriate disjointly supported finite sums thereof.
We consider here gereneralizations of Hindman's Theorem to \emph{arbitrary} semigroups.
Since we do not restrict attention to the abelian case, we usually use multiplicative notation.
Let $S$ be an infinite, finitely colored semigroup. Fix $s\in S$. The homomorphism $n\mapsto s^n$
induces a coloring of $\N$, and by Hindman's Theorem there are
distinct $a_1,a_2,\dots\in\N$ such that all elements
$$s^{a_{i_1}+\dots+a_{i_m}}=s^{a_{i_1}}s^{a_{i_2}}\cdots s^{a_{i_m}}$$
($m\ge 1$, $i_1<i_2<\dots<i_m$) have the same color. Setting $s_n=s^{a_n}$ for all $n$, we have that
all products $s_{i_1}s_{i_2}\cdots s_{i_m}$ ($m\ge 1$, $i_1<i_2<\dots<i_m$) have the same color.
But, unlike Hindman's Theorem, the latter consequence may be trivial: If, for example, $s$ is an idempotent (i.e., $s^2=s$)
then the reason for all products having the same color is that the elements $s_1,s_2,\dots$
and the finite products thereof are all equal to $s$!

Since its publication, several alternative proofs for Hindman's Theorem were published.
The most elegant and powerful one, due to Galvin and Glazer, was first published in Comfort's survey \cite{galvin}.
The Galvin--Glazer proof uses idempotents in the Stone--\v{C}ech compactification $\beta\N$ of $\N$, and generalizes with little effort to a proof of the following theorem.
(Knowledge of the Stone--\v{C}ech compactification is not required in the present paper.)

Say that a semigroup $S$ is \emph{moving} if it is infinite and, for each infinite
$A\sub S$ and each finite $F\sub S$, there are $a_1,\dots,a_k\in A$ such that
$$\{a_1s,a_2s,\dots,a_ks\}\not\sub F$$
for all but finitely many $s\in S$.
Every right cancellative infinite semigroup is moving.
Also, if left multiplication in $S$ is finite-to-one (in particular, if $S$ is left cancellative),
then $S$ is moving.

\bthm[Galvin--Glazer--Hindman]\label{GGH}
Let $S$ be a moving semigroup. For each finite coloring of $S$,
there are \emph{distinct} $a_1,a_2,\dots\in S$ such that all products $a_{i_1}a_{i_2}\cdots a_{i_m}$
($m\ge 1$, $i_1<i_2<\dots<i_m$) have the same color.
\ethm

Our purpose is to generalize the Galvin--Glazer--Hindman Theorem \ref{GGH} to \emph{arbitrary}
infinite semigroups $S$, and to understand the limitations on such generalizations.
We also consider stronger forms of this theorem.

\brem[Attribution]
Theorem \ref{GGH}, which we attribute to Galvin, Glazer and Hindman,
is implicit in Section 4.3 of Hindman and Strauss's monograph \cite{HS}.
There, it is  proved that $S$ is moving if, and only
if, the Stone--\v{C}ech remainder $\beta S\sm S$ is a subsemigroup of $\beta S$.
It follows that $\beta S\sm S$ contains an idempotent, and thus, by the standard
Galvin--Glazer proof of Hindman's Theorem, there are distinct $a_1,a_2,\dots\in S$
as required in Theorem \ref{GGH}.
\erem

\section{Hindman's Theorem everywhere}

As is, the Galvin--Glazer--Hindman Theorem \ref{GGH}
does not generalize to arbitrary semigroups: Consider the following example.

\bexm\label{WhyModFin}
Let $k\in\N$. Let $S$ be the commutative semigroup
$$\{0,1,\dots, k-1\}\cup k\N+1,$$
with the operation of addition modulo $k$.
Assign to each $a\in S$ the color $a \bmod k$.
For all distinct $a_1,a_2,\dots\in S$, we may, by thinning out if necessary,
assume that they are all in $k\N+1$. Consequently, for each $i<k$, $a_1+\dots+a_i=i$, whose color is $i$.
In other words, all colors $i<k$ are obtained when considering all sums of distinct elements from
$\{a_1,a_2,\dots\}$.
\eexm

Thus, we must allow an unbounded finite number of exceptions.
We will soon see that this is the only obstruction to generalizing the
Galvin--Glazer--Hindman Theorem \ref{GGH} to arbitrary semigroups.

We use Shevrin's classification of semigroups.
A semigroup $S$ is \emph{periodic} if $\la s\ra$ is finite for all $s\in S$,
or equivalently, if $\N\not\le S$.
A semigroup $S$ is \emph{right (left) zero} if $ab=b$ ($ab=a$) for all $a,b\in S$.

\bthm[Shevrin \cite{shevrin}]\label{shevrin}
Every infinite semigroup has a subsemigroup of one of the following types:
\be
\itm $(\N,+)$.
\itm An infinite periodic group.
\itm An infinite right zero or left zero semigroup.
\itm $(\N,\vee)$, where $m\vee n:=\max\{m,n\}$.
\itm $(\N,\meet)$, where $m\meet n:=\min\{m,n\}$.
\itm\label{S2} An infinite semigroup $S$ with $S^2$ finite.
\itm The fan semilattice $(\N,\meet)$, with $m\meet n=1$ for distinct $m,n$.
\ee
\ethm

Shevrin's Theorem is stated in \cite{shevrin} in a finer form, replacing \eqref{S2} with
a parameterized list of concrete semigroups.
We will return to this in Section \ref{shevram}.

\bthm\label{HmodFin}
Let $S$ be an infinite semigroup. For each finite coloring of $S$,
there are \emph{distinct} $a_1,a_2,\dots\in S$, and a finite subset $F$ of the (infinite) set
of finite products
$$\FP(a_1,a_2,\dots)=\set{a_{i_1}a_{i_2}\cdots a_{i_m}}{m\ge 1,\, i_1<i_2<\dots<i_m},$$
such that all elements of $\FP(a_1,a_2,\dots)\sm F$ have the same color.
\ethm
\bpf
It suffices to show that every infinite semigroup has a subsemigroup satisfying the assertion of
the theorem.
Apply Shevrin's Theorem \ref{shevrin}. The subsemigroups in cases (1)--(5) are all moving (!),
and thus the Galvin--Glazer--Hindman Theorem \ref{GGH} applies there.

In the remaining cases (6)--(7), let $T$ be the corresponding infinite subsemigroup.
By the pigeon-hole principle, there are distinct $a_1,a_2,\dots\in T$, sharing the same color.
Then
$$\FP(a_1,a_2,\dots) \sub \{a_1,a_2,\dots\}\cup F,$$
where $F$ is $T^2$ in case (6), and $\{1\}$ in case (7),
and thus all elements of $\FP(a_1,a_2,\dots)\sm F$ have the same color.
\epf

\section{Infinite almost-monochromatic subsemigroups}

\bdfn
A colored set $A$ is \emph{monochromatic} if all members of $A$
have the same color.
$A$ is \emph{almost-monochromatic} if all but finitely many members of $A$
have the same color.
\edfn

For which semigroups $S$ is it the case that,
for each finite coloring of $S$, there is an infinite almost-monochromatic subsemigroup of $S$?
We begin with two easy examples.

Let $\Z_2$ be the two element abelian group.
The \emph{direct sum} $\bigoplus_n\Z_2$ is the additive abelian group of all finitely supported elements
of $\Z_2^\N$, with pointwise addition. In other words, $\bigoplus_n\Z_2$ is the group
structure of the countably-infinite-dimensional vector space over the two element field.

\blem\label{lem2b}
For each finite coloring of $\bigoplus_n\Z_2$,
there is an infinite subgroup $H$ of $\bigoplus_n\Z_2$ with $H\sm\{0\}$ monochromatic.
\elem
\bpf
This follows from the Galvin--Glazer--Hindman Theorem \ref{GGH}, since
every group is a moving semigroup, and in the group $\bigoplus_n\Z_2$,
$$\la a_1,a_2,\dots\ra=\set{a_{i_1}+a_{i_2}+\cdots+a_{i_m}}{m\ge 1,\, i_1<i_2<\dots<i_m}\cup\{0\}.\qedhere$$
\epf

\bdfn
A semigroup $S$ is \emph{synchronizing} if $ab\in\{a,b\}$ for all $a,b\in S$.
It is \emph{finitely synchronizing} if there is a finite $F\sub S$ such that
$ab\in\{a,b\}\cup F$ for all $a,b\in S$.
\edfn

Our second example is the class of infinite, finitely synchronizing semigroups.

\blem
Let $S$ be an infinite, finitely synchronizing semigroup.
For each finite coloring of $S$, there is an infinite almost-monochromatic subsemigroup of $S$.
\elem
\bpf
By the pigeon-hole principle, there are distinct $a_1,a_2,\dots\in S$, sharing the same color.
Let $F$ be a finite subset of $S$ such that $ab\in\{a,b\}\cup F$ for all $a,b\in S$.
As
$$\la a_1,a_2,\dots\ra \sub \{a_1,a_2,\dots\}\cup F,$$
$\la a_1,a_2,\dots\ra$ is almost-monochromatic.
\epf

The main result of this section is that the above two easy examples provide a complete answer
to our question. A \emph{$2$-coloring} of a set $A$ is a coloring of $A$ in two colors.

\bthm\label{AlmostMonoSubsemi}
The following are equivalent for semigroups $S$:
\be
\itm For each finite coloring of $S$, there is an infinite almost-mono\-chromatic subsemigroup of $S$.
\itm For each $2$-coloring of $S$, there is an infinite almost-monochro\-matic subsemigroup of $S$.
\itm At least one of the following assertions holds:
\be
\itm $S$ has an infinite, finitely synchronizing subsemigroup.
\itm $\bigoplus_n\Z_2\le S$.
\ee
\ee
\ethm

Item (3)(a) of the theorem may be replaced
by an explicit list of semigroups, namely, the semigroups of types (3)--(7) in Shevrin's Theorem \ref{shevrin}.
Recall that Item (6) can be replaced by a parameterized list of concrete semigroups---see Theorem \ref{shevrin2} below.
Thus, our characterization is completely explicit.

The implication $(3)\Impl (1)$ in Theorem \ref{AlmostMonoSubsemi} is clear. 
Indeed, if $S$ has a subsemigroup $T$ such that for any finite coloring of $T$, 
$T$ contains an infinite almost-monochromatic subsemigroup, then the same holds for $S$. 
The implication $(1)\Impl (2)$  is clear as well. 
The remainder of this section constitutes a proof of the implication $(2)\Impl(3)$.

\blem[Folklore]\label{tech}
Let $G$ be an infinite group such that all elements of $G\sm\{e\}$ have order $2$.
Then $G$ is isomorphic to $\bigoplus_{\alpha\in I} \Z_2$, where $I$ is an index set of
cardinality $|G|$. In particular, $\bigoplus_n \Z_2\le G$.
\elem
\bpf
$G$ is commutative: $[g,h]=ghg^{-1}h^{-1}=(gh)^2=e$ for all $g,h\in G$.
Thus, we may use additive notation for $G$, so that $v+v=0$ for each $v\in G$,
and $G$ is a vector space over the two-element field, necessarily of dimension
$|G|$. In other words, $G$ is isomorphic to  $\bigoplus_{\alpha\in I} \Z_2$.
\epf

\blem\label{Gcolor}
Let $G$ be a group.
There is a $2$-coloring of the elements of $G$ of finite order greater than $2$ such that,
for each coloring of $G$ extending it and each infinite periodic almost-monochromatic
subgroup $H\le G$, $\bigoplus_n\Z_2\le H$.
\elem
\bpf
For each $g\in G$ of finite order greater than $2$, color $g$ and $g^{-1}$ differently.
Let $H$ be an infinite periodic almost-monochromatic subgroup of $G$.
If there are infinitely many $h\in H$ with $h^2\neq e$, then there are infinitely
many such elements of the same color. But then their inverses, which have the opposite color,
also belong to $H$; a contradiction.
Thus, all but finitely many members of $H$ have order $2$.
Let $F$ be the set of elements of order $\neq 2$ in $H$.

Pick $h_1\in H\sm F$. Then $\la h_1\ra=\{h_1,e\}$ is finite.
For $n>1$, assume inductively that all elements of the subgroup
$K:=\la h_1,\dots,h_{n-1}\ra$ of $H$ have order $\le 2$.
Then $K$ is commutative and finite. Pick
$$h_n\in H\sm \bigcup_{h\in K} Fh.$$
Then $h_n\notin K\cup F$ and $h_nh\notin F$ for all $h\in K$.
Consequently, the order of $h_n$ is $2$, and for each $h\in K$, the order of $h_n h$ is $2$.
It follows that $h_nh=hh_n$, and thus $\la h_1,\dots,h_n\ra$ is commutative, finite, and all of its elements have order $\le 2$.

Completing the induction, we have by Lemma \ref{tech}
that $\la h_1,h_2,\dots\ra$ is isomorphic to $\bigoplus_n\Z_2$.
\epf

\blem[Folklore]\label{Ncolor}
There is a $2$-coloring of $\N$ with no infinite almost-monochromatic subsemigroup.
\elem
\bpf
\nc{\g}[1]{\textcolor{green}{#1}}
\renewcommand{\r}[1]{\textcolor{red}{#1}}
Consider the coloring
$$\r1~\g2~\g3~\r4~\r5~\r6~\g7~\g8~\g9~\g{10}~\r{11}~\r{13}~\r{14}~\r{15}~\r{16}~
\g{17}~\g{18}~\g{19}~\g{20}~\g{21}~\g{22}~\r{23}~\dots,$$
where the lengths of the intervals of elements of identical colors are $1,2,3,\dots$.
For each $n\in\N$, $\la n\ra$ intersects every monochromatic interval of length $\ge n$.
\epf

For a semigroup $S$ and an idempotent $e\in S$, 
let $G(e)$ be the maximal subgroup of the semigroup $S$ containing the idempotent $e$.
As groups have exactly one idempotent,
$G(e_1)\cap G(e_2)=\emptyset$ for all distinct idempotents $e_1,e_2\in S$.

\theoremstyle{plain}
\newtheorem{tclem}[thm]{True Color Lemma}

\begin{tclem}\label{truecolor}
For each semigroup $S$, there is a $2$-coloring of $S$ such that:
\be
\itm Every almost-monochromatic subsemigroup of $S$ is periodic; and
\itm Every infinite almost-monochromatic subgroup of $S$ contains $\bigoplus_n\Z_2$
as a subgroup.
\ee
\end{tclem}
\bpf
An \emph{orbit} in $S$ is a subset of the form $\la s\ra$ for some $s\in S$.
If there are infinite orbits in $S$, use Zorn's Lemma to obtain a maximal family $\cF$ of disjoint infinite orbits in $S$.
If there are none, let $\cF=\emptyset$.
For each $\la s\ra\in\cF$, $\la s\ra$ is isomorphic to $\N$.
Use Lemma \ref{Ncolor} to obtain,
for each $\la s\ra\in\cF$, a coloring of $\la s\ra$ in red and green, such that
$\la s\ra$ has no almost-monochromatic subsemigroup.

Let $e$ be an idempotent of $S$. The elements of finite order greater than $2$ in $G(e)$
do not belong to an infinite orbit, and are thus not colored yet. Color these elements in red
and green, as in Lemma \ref{Gcolor}.
As the groups $G(e)$ are disjoint for distinct idempotents, this can be done for
all idempotents.

Extend our partial $2$-coloring to an arbitrary $2$-coloring of $S$.

(1) Let $T$ be a non-periodic subsemigroup of $S$. Pick $t\in T$ with $\la t\ra$ infinite.
By the maximality of $\cF$, $\la t\ra$ intersects some $\la s\ra\in\cF$. Let $n$ be such that
$t^n\in \la s\ra$. Then the subsemigroup $\la t^n\ra$ of $\la s\ra$ is not almost-monochromatic.
In particular, $T$ is not almost-monochromatic.

(2) Let $G$ be an infinite almost-monochromatic subgroup of $S$. By (1), $G$ is periodic.
Let $e$ be the idempotent of $G$. Then $G\le G(e)$, and by Lemma
\ref{Gcolor}, $\bigoplus_n\Z_2\le G$.
\epf

\bpf[Proof of Theorem \ref{AlmostMonoSubsemi}]
Assume that, for each $2$-coloring of $S$,
there is an infinite almost-monochromatic subsemigroup of $S$.
Color $S$ as in the True Color Lemma \ref{truecolor}.
Let $T$ be an almost-monochromatic subsemigroup of $S$.
By the True Color Lemma, $T$ is periodic.

If $T$ has no infinite subgroup, then cases (1) and (2) in Shevrin's Theorem \ref{shevrin} are excluded.
As each of the semigroups in the remaining cases of Shevrin's Theorem is finitely synchronizing,
$S$ has an infinite, finitely synchronizing subsemigroup.

And if $T$ has an infinite subgroup, $G$, then by the True Color Lemma,
$\bigoplus_n\Z_2\le G$.
\epf

The case of \emph{groups} is of independent interest.

\bthm\label{2b}
The following are equivalent for groups $G$:
\be
\itm For each finite coloring of $G$, there is an infinite almost-mono\-chromatic subgroup of $G$.
\itm For each finite coloring of $G$, there is an infinite almost-mono\-chromatic subsemigroup of $G$.
\itm For each $2$-coloring of $G$, there is an infinite almost-monochro\-matic subgroup of $G$.
\itm For each $2$-coloring of $G$, there is an infinite almost-monochro\-matic subsemigroup of $G$.
\itm $\bigoplus_n\Z_2\le G$.
\ee
\ethm
\bpf
Clearly, the implications $(1\Impl 2)$, $(1\Impl 3)$, $(2\Impl 4)$ and $(3\Impl 4)$ hold.

$(4\Impl 5)$ Apply Theorem \ref{AlmostMonoSubsemi} to the group $G$.
If $T$ is an infinite, finitely synchronizing subsemigroup of $G$, then $T$ is a periodic
subsemigroup of a group, and thus a group. But infinite groups cannot be finitely synchronizing.
Indeed, let $F$ be a finite subset of $G$, $a\in G\sm \{e\}$.
Since left multiplication by $a$ is injective, there is $b\in G\sm\{e\}$ such that
$ab\notin F$. As $a,b\neq e$, $ab\notin\{a,b\}$, and thus $a,b\notin\{a,b\}\cup F$.
Consequently, we are in case (3.b) of Theorem  \ref{AlmostMonoSubsemi}, that is, $\bigoplus_n\Z_2\le G$.

$(5\Impl 1)$ Lemma \ref{lem2b}.
\epf

\section{Unordered products}

For distinct $a_1,a_2,\dots\in S$, let
$$\FPhat(a_1,a_2,\dots)=\set{a_{i_1}a_{i_2}\cdots a_{i_m}}{m\ge 1,\, i_1,i_2,\dots,i_m
\mbox{ are distinct}}.$$
We apply the information gathered in the previous sections
to the following question: Let $S$ be a prescribed infinite semigroup.
Is it true that, for each finite coloring of $S$, there are distinct $a_1,a_2,\dots\in S$ such that $\FPhat(a_1,a_2,\dots)$ is almost-monochromatic?

To see that the new question is different than the one studied in the previous section,
note that, by Hindman's Theorem, for each finite coloring of $\N$, there are distinct $a_1,a_2,\dots\in\N$ such that
the set
$$\FPhat(a_1,a_2,\dots)=\set{a_{i_1}+a_{i_2}+\cdots+a_{i_m}}{m\ge 1,\, i_1,i_2,\dots,i_m
\mbox{ are distinct}}$$
is monochromatic, but there is a $2$-coloring of $\N$ with no infinite
almost-monochromatic subsemigroup (Lemma \ref{Ncolor}).

\bthm\label{FPhatThm}
Assume that $S$ has an infinite subsemigroup 
with no infinite, finitely generated, periodic subgroup.
For each finite coloring of $S$, there are distinct $a_1,a_2,\dots\in S$ with
$\FPhat(a_1,a_2,\dots)$ almost-monochromatic.
\ethm
\bpf
Assume that the theorem fails for $S$. Then $(\N,+)$ is not a subsemigroup of $S$. 
By moving to a subsemigroup of $S$, if needed,
we may assume that $S$ is an infinite semigroup 
with no infinite, finitely generated, periodic subgroup.

By Theorem \ref{AlmostMonoSubsemi}, $S$ does not contain an infinite
finitely synchronizing subsemigroup.
Thus, by Shevrin's Theorem \ref{shevrin}, $S$ has an infinite periodic subgroup $G$.

If $G$ is locally finite, then it contains an infinite \emph{abelian} group $H$ \cite{hall}.
As groups are moving, by the Galvin--Glazer--Hindman Theorem \ref{GGH}
there are distinct $a_1,a_2,\dots\in H$ such that $\FPhat(a_1,a_2,\dots)$ is monochromatic;
a contradiction.

Thus, $G$ is not locally finite. Let $F\sub G$ be a finite set with $H:=\la F\ra$ infinite.
Then $H$ is an infinite, finitely generated, periodic subgroup of $S$; a contradiction.
\epf

The condition on $S$ in Theorem \ref{FPhatThm} is quite mild:
The 1902 \emph{Burnside Problem} \cite{burn}, that asked whether
there is, at all, an infinite finitely generated periodic group,
was only answered (in the affirmative) in 1964 \cite{golod}.

The question whether the condition in Theorem \ref{FPhatThm} can be eliminated
is equivalent to a 1978 problem of Milliken.

\bprp
Let $G$ be an infinite, finitely colored group.
Assume that $a_1,a_2,\dots\in G$ are distinct elements such that $\FPhat(a_1,a_2,\dots)$ is almost-monochromatic.
Then there is a subsequence $a_{i_1},a_{i_2},\dots$ of $a_1,a_2,\dots$ such that
$\FPhat(a_{i_1},a_{i_2},\dots)$ is monochromatic.
\eprp
\bpf
Let $F$ be the finite set of elements of $\FPhat(a_1,a_2,\dots)$ having exceptional colors.
Pick $a_{i_1}\in\{a_1,a_2,\dots\}\sm F$, and for $n>1$, let $P$ be the set of all products of
at most $n-1$ distinct elements from $\{a_{i_1},\dots,a_{i_{n-1}}\}$, including also $e$.
Pick $a_{i_n}\in \{a_1,a_2,\dots\}\sm P^{-1}FP^{-1}$, with $i_n>i_{n-1}$.
\epf

\bprb[Milliken \cite{milliken}]
Is it true that, for each infinite, finitely colored group $G$, there are distinct $a_1,a_2,\dots\in G$
such that $\FPhat(a_1,a_2,\dots)$ is monochromatic?
\eprb

In 1968, Novikov and Adian \cite{NA1} proved that, for each $m \ge 2$ and each large enough 
odd $n$, the \emph{Burnside group}
$$G=\la x_1,\dots,x_m\,:\,x^n=1\ra$$
(where $x^n=1$ for all $x\in G$) is infinite (cf.\ Adian \cite{a2}).
As was already noted by Milliken \cite{milliken},
for large enough odd $n$ these groups
have no infinite abelian subgroups \cite{NA2}, and thus the Galvin--Glazer--Hindman
Theorem does not apply to them directly.

A group $G$ is a \emph{Tarski Monster} if, for some prime number $p$, all
proper subgroups of $G$ have cardinality $p$.
Tarski Monsters exist for all large enough primes $p$ (Olshanskii \cite{tarski}; cf.\
Adian--Lys\"enok \cite{al}).
Clearly, Tarsky Monsters do not have infinite abelian subgroups.
Thus, it may be possible to address Milliken's problem by finding the ``true color'' of
some Tarski Monster\dots

\section{A semigroup structure theorem of Shevrin, via Ramsey's Theorem}\label{shevram}

In the previous sections, we applied Shevrin's theory to coloring theory. We conclude with
an application in the converse direction.

The following assertion is made in \cite{shevrin}.
For completeness, we give a proof.

\blem[Shevrin \cite{shevrin}]\label{she}
Let $S$ be a semigroup generated by $A$, such that, for some natural numbers $h>1$ and $d$:
\be
\itm[(a)] $abc=def$ for all $a,b,c,d,e,f\in A$;
\itm[(b)] $a^h=a^{h+d}$ for all $a\in A$.
\ee
Then:
\be
\itm $S^3={\la a \ra}^3$ for each $a\in A$.
\itm $S^3$ is finite.
\itm There is a unique idempotent $e\in S$.
\itm For all $a,b\in A$, $ae=be$.
\ee
\elem

\bpf
(1) Each $s\in S^3$ is a product of $k\ge 3$ elements of $A$. Applying (a) repeatedly,
we conclude that $s=a^k$.

(2) Fix $a\in A$. By (b), $\la a\ra $ is finite. Apply (1).

(3) Fix $a\in A$. By (b), $\la a\ra$ is finite, and thus there is an idempotent $e=a^k$ in $\la a\ra$.
Let $s\in S$ be an idempotent. By (1), $s=s^3\in\la a\ra$, say $s=a^m$. Then $s=s^k=a^{mk}=e^m=e$.

(4) Let $a,b\in A$. By (1), $e=e^3\in \la a\ra^3$, and hence $e=a^k$ for some $k\geq 3$.
By (a), $ae=a^{k+1}=a^3a^{k-2}=ba^2a^{k-2}=ba^k=be$.
\epf

Following Shevrin \cite{shevrin}, say that a semigroup $S$ is of
\emph{type $[h,d]$} for $h,d\in\N$ with $h>1$,
if $S$ is generated by a countably infinite alphabet $x_1,x_2,\dots$, with the following defining relations:
\be
\itm[(HD1)] ${x_i}^2={x_1}^2$ for all $i$;
\itm[(HD2)] $x_{i}x_{j}=x_{1}x_{2}$ and $x_{j}x_{i}=x_{2}x_{1}$ for all $i<j$;
\itm[(HD3)] $x_{i}x_{j}x_{k}={x_1}^3$ for all $i,j,k$;
\itm[(HD4)] ${x_i}^h={x_i}^{h+d}$ for all $i$;
\ee
and possibly by additional relations, equating some or all of the elements:
${x_1}^2$, $x_1x_2$, $x_2x_1$, $(x_1e)^2$,
where $e$ is the unique idempotent of $S$ (Lemma \ref{she}(3)).

Shevrin proves in \cite{shevrin} a finer version of Theorem \ref{shevrin},
where ``An infinite semigroup $S$ with $S^2$ finite'' is replaced by ``A semigroup of type $[h,d]$.''
In the course of his proof, however, he essentially proves the equivalence of these two versions of
Theorem \ref{shevrin}.
We give a short, complete proof using Ramsey's celebrated coloring theorem. Ramsey's Theorem
asserts that, for each finite coloring of the edges of an infinite complete graph,
there is an infinite complete subgraph with all edges of the same color.

We first treat the easier implication of Shevrin's Theorem.

\bprp[Shevrin \cite{shevrin}]
Let $S$ be a semigroup of type $[h,d]$. Then $S$ is infinite, and $S^2$ is finite.
\eprp
\bpf
As all of the words in the defining relations have more than one letter,
there are no relations applicable to a single letter. Consequently, all letters of $A$ are
distinct in $S$, and $S$ is infinite.

By Lemma \ref{she}(2), $S^3$ is finite. By the defining relations,
$S^2\sm S^3\sub \{{x_1}^2$, $x_1x_2$, $x_2x_1\}$. Thus, $S^2$ is finite.
\epf

\bthm[Shevrin \cite{shevrin}]\label{shevrin2}
Let $S$ be an infinite semigroup with $S^2$ finite.
Then $S$ has a subsemigroup of type $[h,d]$, for some natural numbers $h>1$ and $d$.
\ethm
\bpf
As $S^3\sub S^2$, $S^3$ is finite too.
Pick distinct elements $a_1,a_2,\dots\in S\sm S^2$.
Consider the complete infinite graph with vertex set $V=\{a_1,a_2,\dots\}$.
Think of the finite set $S^3\times S^2\times S^2\times S^2$ as a set of colors,
and define a finite coloring of the edges of our graph,
$$c:[V]^2\rightarrow S^3\times S^2\times S^2\times S^2,$$
by
$$c(\{a_i,a_j\}):=({a_i}^3,{a_i}^2,a_ia_j,a_ja_i).$$
for all $i<j$.
By Ramsey's Theorem, there are $i_1<i_2<\dots$ such that all edges among the vertices in the set
$\{a_{i_1},a_{i_2},\dots\}$ have the same color. Denote $b_n=a_{i_n}$ for all $n$.

Let $1\le i<j$. Then
$$({b_{i}}^3,{b_{i}}^2,b_{i}b_{j},b_{j}b_{i})=c(\{b_{i},b_{j}\})=c(\{b_2,b_{3}\})=
({b_2}^3,{b_2}^2,b_2b_{3},b_{3}b_2).$$
Hence, for all $1\le i<j$,
\begin{eqnarray*}
{b_i}^3 & = & {b_2}^3;\\
{b_i}^2 & = & {b_2}^2;\\
b_{i}b_{j} & = & b_2b_3;  \mbox{ and}\\
b_jb_i & = & b_3b_2.
\end{eqnarray*}
We claim that the subsemigroup $T=\la b_2,b_3,\dots \ra$ is of type $[h,d]$ for some $h>1$ and $d$,
with respect to the alphabet $b_2,b_3,\dots$.

The elements $b_2,b_3,\dots$ are distinct, being a subsequence of the sequence $a_1,a_2,\dots$
of distinct elements.

We have already proved that Relations (HD1) and (HD2) hold.

(HD3) Fix $k\ge 2$. If $i<j$, then $b_ib_j=b_2b_3=b_1b_k$, and thus
$$b_ib_jb_k=b_1b_kb_k=b_1{b_2}^2={b_1}^3={b_2}^3.$$
If $i>j$, then
\begin{eqnarray*}
b_ib_jb_k & = & (b_3b_2)b_k=(b_{k+1}b_k)b_k = b_{k+1}{b_k}^2=\\
& = & b_{k+1}{b_2}^2=b_{k+1}{b_{k+1}}^2={b_{k+1}}^3={b_2}^3.
\end{eqnarray*}
If $i=j$, then
$$b_ib_jb_k=(b_ib_i)b_k={b_2}^2b_k={b_{k}}^2b_k={b_{k}}^3={b_2}^3.$$

(HD4) Denote $b=b_2$. As $S^2$ is finite, so is $\la b\ra$. Take minimal $h$ and $d$ such that
$b^h=b^{h+d}$. As $b\in S\sm S^2$, $h>1$. Thus, for all $i\ge 2$, we have by (HD1) and (HD3) that
$${b_i}^h=b^h=b^{h+d}={b_i}^{h+d}.$$

By (HD1)--(HD4),
$$T=\{b_2,b_3,\dots\}\cup\{{b_2}^2, b_2b_3, b_3b_2\}\cup\{b^3,\dots, b^{h+d-1}\},$$
where $b=b_2$. (In the case $h=2$ and $d=1$, the rightmost set in this union is empty.)

By Lemma \ref{she}, (HD3) and (HD4), there is a unique idempotent $e$ in $T$.
It thus remains to show that no additional equalities, except perhaps ones among
${b_2}^2$, $b_2b_3$, $b_3b_2$, and $(be)^2$, hold.

We already observed that the elements $b_2,b_3,\dots$ are distinct, and by their choice,
do not belong to $S^2$. Thus, equalities may only hold among members of the set
$\{{b_2}^2, b_2b_3, b_3b_2\}\cup\{b^3,b^4,\dots,b^{h+d-1}\}$.
In the case $h=2$ and $d=1$, the rightmost set is empty, and we are done.
Consider the other cases.
By the minimality of $h$ and $d$, the elements $b=b_2,b^2={b_2}^2,b^3,\dots,b^{h+d-1}$ are distinct,
and
$$G:=\{b^h,b^{h+1},\dots,b^{h+d-1}\}$$
is a group.
The idempotent element of $G$ must be the unique idempotent $e$ of $S$.

Thus, equalities may only hold among $b_2b_3,b_3b_2$ (which is fine), or be
of the form $b^n=b_2b_3$ or $b^m=b_3b_2$,  for some (necessarily unique)
$2\le n,m\le h+d-1$. It suffices to show that the former case is equivalent to $b_2b_3={b_2}^2$ or to $b_2b_3=(be)^2$ and the latter to $b_3b_2={b_2}^2$ or to $b_3b_2=(be)^2$.

We prove the assertion for $b_2b_3$; the other proof being identical.

As $e\in G$, $be\in G$ and thus so is $(be)^2$. Thus, if $b_2b_3=(be)^2$, then there is $2\le h\le n\le h+d-1$ with $b_2b_3=b^n$.

For the other direction, if $b^n=b_2b_3$ for $n>2$ then
$$b^{n+1}=b^nb=b_2b_3b=b^3,$$
and therefore $h\le 3\le n$. Thus, $b^n\in G$. As $e\in G\sub \la b \ra$, $be=eb$. Then
$$b^n(be)=b^{n+1}e=b^3e=b^3e^3=(be)^2(be).$$
As $b^n,be\in G$, this implies that $b_2b_3=b^n=(be)^2$, as required.
\epf

\subsection*{Acknowledgments}
We thank the referee for his or her suggestions, that made the paper more accessible.

\ed